\documentclass[12pt]{article}
\usepackage{amssymb}
\usepackage{amsmath}
\usepackage{amsbsy}
\usepackage{amsthm}
\usepackage{rotating}

\newcommand{\be}{\begin{equation}}
\newcommand{\ee}{\end{equation}}
\newcommand{\ba}{\begin{eqnarray}}
\newcommand{\ea}{\end{eqnarray}}
\newcommand{\ban}{\begin{eqnarray*}}
\newcommand{\ean}{\end{eqnarray*}}

\newtheorem{theo}{Theorem}[section]
\newtheorem{prop}[theo]{Proposition}

\newtheorem{cor}[theo]{Corollary}
\newtheorem{defi}[theo]{Definition}

\newtheorem{conj}{Conjecture}[section]
\begin{document}

\title{The Novikov conjecture for algebraic K-theory of the group algebra over the ring of Schatten class operators}
\author{Guoliang Yu\footnote{The author is  partially supported by NSF.}}
\date{ }

\maketitle

{\noindent  Abstract:} In this paper, we prove the algebraic K-theory Novikov conjecture for group algebras
over the ring of Schatten class operators. The main technical tool in the proof is an explicit construction of the Connes-Chern character.

\section{Introduction}

Let  $\Gamma$ be  a group and $R$ be an H-unital ring. Let $R\Gamma$ be the group algebra of the group $\Gamma$ over the ring $R$.
 The isomorphism conjecture of Farrell-Jones  states that the following assembly map is
an isomorphism:

$$A: H_n^{Or\Gamma}(E_{\cal VCY} (\Gamma), {\mathbb K}(R)^{-\infty}) \longrightarrow K_n (R\Gamma),$$
where ${\cal VCY}$ is the family of virtually cyclic subgroups of $\Gamma$,  $ E_{\cal VCY} (\Gamma)$ is the universal $\Gamma$-space with
isotropy in ${\cal VCY}$, $ H_n^{Or\Gamma}(E_{\cal VCY} (\Gamma), {\mathbb K}(R)^{-\infty})$ is a generalized $\Gamma$-equivariant homology theory
associated to the non-connective algebraic K-theory spectrum $ {\mathbb K}(R)^{-\infty}$, and $K_n (R\Gamma)$ is the algebraic K-theory of $R\Gamma$.

The isomorphism conjecture provides an algorithm for computing the algebraic K-theory of $R\Gamma$ in terms of the algebraic K-theory of $R$. This  conjecture was introduced in [FJ1] for $R={\mathbb Z}$ and for unital rings $R$ in [BFJR]. When $R$ is H-unital, the isomorphism conjecture follows from the unital case by using the excision theorem in algebraic K-theory  [SW]. The algebraic K-theory isomorphism  conjecture goes back to [H]. There are analogous conjectures in L-theory [Q1] [Q2] and $C^\ast$-algebra K-theory [BC1]. Important cases of the isomorphism conjecture have been verified in [FJ1] [FJ2] and [BLR].

The algebraic K-theoretic Novikov conjecture states  that the assembly map:
$$ H_n (B\Gamma, {\mathbb K}(R)^{-\infty}) \longrightarrow K_n (R\Gamma),$$
is rationally injective, where $B\Gamma$ is the classifying space of the group $\Gamma$. The algebraic K-theoretic Novikov conjecture follows from the (rational) injectivity part of the isomorphism conjecture.
By a remarkable theorem of B\"{o}kstedt-Hsiang-Madsen [BHM], the algebraic K-theoretic Novikov conjecture holds for  $R=\mathbb{Z}$ if the  homology groups
of $\Gamma$ are finitely generated.

The main purpose of this paper is to prove the (rational) injectivity part of the algebraic K-theory isomorphism conjecture for group algebras
over the ring of Schatten class operators. As a consequence, we obtain the algebraic K-theory Novikov  conjecture for group algebras
over the ring of Schatten class operators. The motivation for considering  group algebras over the ring of Schatten class operators comes from the deep work of Connes-Moscovici on higher index theory of elliptic operators and its applications to the Novikov conjecture [CM]. In Connes-Moscovici's higher index theory,  the K-theory of the group algebra over the ring of Schatten class operators  serves as  the receptacle for the higher index of an elliptic operator.

For the convenience of readers we recall that, for any $p\geq 1$, an operator $T$ on an infinite dimensional and separable Hilbert space $H$ is said to be Schatten $p$-class if $tr ((T^\ast T)^{p/2} )<\infty$, where $tr$ is the standard trace defined by $tr (P)=\sum_{n} <Pe_n, e_n>$ for any bounded operator
 $P$ acting on $H$ and an orthonormal basis $\{e_n\}_n$ of $H$  ($tr (P)$ is independent of the choice of the orthonormal basis).
Let ${\cal S}_p$ be the ring of all  Schatten $p$-class operators on an infinite dimensional and separable Hilbert space. We define the ring ${\cal S}$ of all Schatten class operators to be  $\cup_{p\geq 1} {\cal S}_p$.

\pagebreak

The  following theorem is the main result of this paper.

\begin{theo} Let ${\cal S}$ be the ring of all  Schatten class operators on an infinite dimensional and separable Hilbert space. The assembly map
$$A: H_n^{Or\Gamma}(E_{\cal VCY} (\Gamma), {\mathbb K}({\cal S})^{-\infty}) \longrightarrow K_n ({\cal S}\Gamma)$$ is rationally injective
for any   group $\Gamma$,
where ${\cal S}\Gamma$ is the group algebra of the group $\Gamma$ over the ring ${\cal S}$.
\end{theo}

As a consequence, we obtain the algebraic K-theoretic Novikov conjecture for the group algebra ${\cal S} \Gamma$.

The main technical tool in the proof of Theorem 1.1 is an explicit construction of a Connes-Chern character using  an equivariant cyclic simplicial homology theory. As a consequence of this explicit construction, we obtain a local property  of the Connes-Chern character. This local property of the Connes-Chern character plays an important role in the proof.

This paper is organized as follows. In Section 2, we collect a few preliminary results which will be used later in the paper.
In Section 3, we reduce our main theorem to the case of lower algebraic K-theory. In Section 4, we introduce a cyclic simplicial
 homology theory to construct a Connes-Chern character. The Connes-Chern character plays a crucial role in the proof of the main theorem.
 We use the explicit construction of the Connes-Chern character to prove an important local property of the Connes-Chern character for K-theory elements with small propagation.
In Section 5, we prove the main theorem of this paper.

The author wishes to thank Alain Connes, Max Karoubi, Xiang Tang, Andreas Thom, Shmuel Weinberger,  and  Rufus Willett for   inspiring discussions  and very helpful comments.  In particular, the author would like to express his gratitude to Guillermo  Corti\~nas for his  detailed comments about the paper and several stimulating discussions. We would like to mention that Guillermo  Corti\~nas  and Giesela Tartaglia have given a new proof of Theorem 1.1 in [CTa]. Part of this work was done at Shanghai Center for Mathematical Sciences (SCMS)  and the author wishes to thank SCMS for providing excellent working environment.

\section{Preliminaries}

In this section, we collect a few   concepts and results useful for this paper.

Let $R$ be a ring and let $R^{+}$ be the unital ring obtained from $R$ by adjoining a unit.
The ring $R$ is defined to be H-unital if $Tor_i^{R^{+}} ({\mathbb Z}, {\mathbb Z})=0$ for all $i$.
The importance of H-unitality is that it guarantees excision in algebraic K-theory [SW].

If $R$ is a ${\mathbb Q}$-algebra and
$R^{+}_{{\mathbb Q}  }$ is the unital ${\mathbb Q}$-algebra obtained from $R$ with the unit adjoined, then
 $R$ is H-unital if and only if $Tor_i^{R^{+}_{\mathbb Q}} ({\mathbb Q}, {\mathbb Q})=0$ [SW].

The following result follows from [W1] and Theorem 8.2.1 of [CT].

\begin{theo} ${\cal S}$ is H-unital.
\end{theo}

By Theorem 7.10 in [SW], we have:

\begin{theo} If $R$ is H-unital, then $R\Gamma$ is H-unital for any group $\Gamma$.
\end{theo}

As a consequence, we obtain that ${\cal S}\Gamma$ is H-unital.

Recall that a ring $R$ is called $K_n$-regular if the natural map: $$K_n(R)\rightarrow K_n( R[t_1, \cdots, t_m]),$$
 is an isomorphism for each $m\geq 1$. We say that $R$ is $K$-regular if $R$ is $K_n$-regular for all $n$.

The following result is a special case of Theorem 8.2.5 in [CT].

\begin{theo} ${\cal S}$ is $K$-regular.
\end{theo}

The following result follows from the proof of Proposition 2.14 in [LR].

\begin{prop}
If $R$ is a $K$-regular ${\mathbb R}$-algebra, then the natural map:
$$ H_n^{Or\Gamma}(E_{\cal FIN} (\Gamma), {\mathbb K}( R)^{-\infty}) \rightarrow
 H_n^{Or\Gamma}(E_{\cal VCY} (\Gamma), {\mathbb K}(R)^{-\infty}),$$
is an isomorphism, where
${\cal FIN}$ is  the family
of finite subgroups of $\Gamma$ and $ E_{\cal FIN} (\Gamma) $ is the universal $\Gamma$-space with isotropy in ${\cal FIN}$.
\end{prop}

The above proposition implies that the isomorphism conjecture for the ring ${\cal S}$ is equivalent to the statement that
the assembly map:
$$A: H_n^{Or\Gamma}(E_{\cal FIN} (\Gamma), {\mathbb K}({\cal S})^{-\infty}) \rightarrow K_n({\cal S}\Gamma),$$
is an isomorphism and Theorem 1.1 is equivalent to the statement that the above assembly map is rationally injective.

By a result in [CT], we know that $K_n({\cal S})$ is 2-periodic and $K_0 ({\cal S})={\mathbb Z}$ and $K_1 ({\cal S}) =0$. This implies that the domain of the assembly map $A$ in Theorem 1.1  is rationally isomorphic to $\oplus _{k~even} H_{n+k}^{Or\Gamma}(E_{\cal FIN} (\Gamma), {\mathbb Q}).$

\section{Reduction to the lower algebraic $K$-theory case}

In this section, we prove the following reduction result.

\begin{prop} Theorem 1.1 follows from the following special case of  the theorem  for lower algebraic  K-theory: i.e.
the assembly map
$$A: H_n^{Or\Gamma}(E_{\cal VCY} (\Gamma), {\mathbb K}({\cal S})^{-\infty}) \rightarrow K_n ({\cal S}\Gamma)$$ is rationally injective
for $n\leq 0$.
\end{prop}
\proof By Proposition 7.2.3, Remark 7.2.6 and Theorem 8.2.5 in [CT], we have $K_n ({\cal S}) ={\mathbb Z}$ when $n$ is even and $K_n ({\cal S}) =0$ when $n$ is odd.
It follows that
$$H_{-2} ({pt}, {\mathbb K}({\cal S})^{-\infty})={\mathbb Z}.$$
By definition, the assembly map:
$$ A: H_{-2} ({pt}, {\mathbb K}({\cal S})^{-\infty})\rightarrow K_{-2}({\cal S})$$
is  an isomorphism  and it maps the generator $z$ of $H_{-2} ({pt}, {\mathbb K}({\cal S})^{-\infty})$ to the Bott element of $ K_{-2}({\cal S})$
(denoted by $b$). For any positive integer $k$, we can use the product operation to construct the Bott element $b^k$ in $K_{-2k}({\cal S}),$
where the product is defined using a natural (injective) homomorphism $ {\cal S}\otimes {\cal S}\rightarrow {\cal S}$ induced by  the homomorphism $ {\cal S(H)}\otimes {\cal S(H)}\rightarrow {\cal S(H\otimes H)}$ and  the isomorphism ${\cal S(H\otimes H)}\cong {\cal S(H)}$ (here ${\cal H}$ is an infinite dimensional and separable Hilbert space, ${\cal S(H)}$ and ${\cal S(H\otimes H)}$ respectively denote the rings of Schatten class operators on ${\cal H}$ and ${\cal H\otimes H}$,  and  $ {\cal S}\otimes {\cal S}$ is the algebraic tensor product of $ {\cal S}$ with $ {\cal S}$).

When $n=2k$, we have the following commutative diagram:

$$\begin{array}[c]{ccc}
H_n ^{Or \Gamma}(E_{\cal VCY} (\Gamma), {\mathbb K}({\cal S})^{-\infty}) &\stackrel{A}{\rightarrow}& K_{n}({\cal S}\Gamma)\\
\downarrow\scriptstyle{\times z^k}&&\downarrow\scriptstyle{\times b^k}\\
H_0 ^{Or \Gamma}(E_{\cal VCY} (\Gamma), {\mathbb K}({\cal S})^{-\infty})&\stackrel{A}{\rightarrow}& K_0 ({\cal S}\Gamma)
\end{array},$$
where the vertical product maps are well defined with the help of a natural homomorphism $ {\cal S}\otimes {\cal S}\rightarrow {\cal S}$ and the $K$-theory properties of ${\cal S} $ and ${\cal S} \Gamma$  being $H$-unital (Theorems 2.1 and 2.2).
 By Theorems 8.2.5 and 6.5.3 in [CT] and Theorem 8.3 (the Bott periodicity)   in [Cu2], we know that the Bott element $b^k$ is a generator of $K_{-2k}({\cal S}).$
   It follows that  the  product map
 $$ H_i ( {pt}, {\mathbb K}({\cal S})^{-\infty}) \stackrel{ \times z^k}{\longrightarrow} H_{i-2k} ( {pt}, {\mathbb K}({\cal S})^{-\infty})$$
 is an isomorphism for every integer $i$.  This implies that the  product map
 $$H_i ^{Or \Gamma}(E_{\cal VCY} (\Gamma), {\mathbb K}({\cal S})^{-\infty}) \stackrel{ \times z^k}{\longrightarrow} H_{i-2k} ^{Or \Gamma}(E_{\cal VCY} (\Gamma), {\mathbb K}({\cal S})^{-\infty})$$ is an isomorphism by using the fact that
   both homology theories \linebreak $\{ H_i ^{Or \Gamma}(~\cdot~, {\mathbb K}({\cal S})^{-\infty})\}_{i\in {\mathbb Z}}$ and $\{ H_{i-2k} ^{Or \Gamma}(~\cdot~, {\mathbb K}({\cal S})^{-\infty})\}_{i\in {\mathbb Z}}$ have a   Mayer-Vietoris sequence and a five lemma argument.

When $n=2k+1$, we have the following commutative diagram:

$$\begin{array}[c]{ccc}
H_n ^{Or \Gamma}(E_{\cal VCY} (\Gamma), {\mathbb K}({\cal S})^{-\infty}) &\stackrel{A}{\rightarrow}& K_{n}({\cal S}\Gamma)\\
\downarrow\scriptstyle{\times z^{k+1}}&&\downarrow\scriptstyle{\times b^{k+1}}\\
H_{-1} ^{Or \Gamma}(E_{\cal VCY} (\Gamma), {\mathbb K}({\cal S})^{-\infty})&\stackrel{A}{\rightarrow}& K_{-1} ({\cal S}\Gamma)
\end{array},$$
where  the vertical product maps are well defined  with the help of a natural homomorphism $ {\cal S}\otimes {\cal S}\rightarrow {\cal S}$ and
the $K$-theory properties of ${\cal S} $ and ${\cal S} \Gamma$  being $H$-unital (Theorems 2.1 and 2.2).
 By the same argument as in the even case, we know that  the  product map
 $$ H_i ( {pt}, {\mathbb K}({\cal S})^{-\infty}) \stackrel{ ~~\times z^{k+1}}{\longrightarrow} H_{i-(2k+2)} ( {pt}, {\mathbb K}({\cal S})^{-\infty})$$
 is an isomorphism for every integer $i$.
  This fact, together with  a standard Mayer-Vietoris sequence and five lemma argument, implies that the product map
   $$H_i ^{Or \Gamma}(E_{\cal VCY} (\Gamma), {\mathbb K}({\cal S})^{-\infty}) \stackrel{~~ \times z^{k+1}}{\longrightarrow} H_{i-(2k+2)} ^{Or \Gamma}(E_{\cal VCY} (\Gamma), {\mathbb K}({\cal S})^{-\infty})$$ is an isomorphism.

Now Proposition 3.1 follows from the above commutative diagrams and the fact that the left vertical maps in the diagrams are isomorphisms.
\qed

\section{Cyclic simplicial homology theory and the  Connes-Chern character}

In this section, we introduce an equivariant cyclic simplicial homology theory to
construct the Connes-Chern character for  $K_n ({\cal S}\Gamma)$ when $n\leq 0$. The Connes-Chern character is a key tool in the proof of the main theorem.  We use this explicit construction to prove an important local property of the Connes-Chern character for K-theory elements with small propagation. This local property will be useful in the proof of the main theorem.

Let $X$ be a simplicial complex. Let $\sigma$ be a simplex of $X$. Define two orderings of its vertex set to be equivalent if they differ
from each other by an even permutation. Each of the equivalence classes is called an orientation of $\sigma.$
 If $\{v_0, ..., v_k\}$ is the set of all vertices of $\sigma$, we use the symbol
$[v_0, \cdots, v_k]$ to denote the oriented simplex with the particular ordering $(v_0,\cdots, v_k)$.

A locally finite $k$-chain on $X$ is a formal sum
$$\sum_{(v_0, \cdots, v_k)} c_{(v_0, \cdots, v_k)} [v_0, \cdots, v_k],$$
where

\noindent (1) the summation is taken over all  orderings $ (v_0, \cdots, v_k)$ of  all $k$-simplices $\{ v_0, \cdots, v_k \}$ of   $X$ and
%${\mathbb R}$-valued  function $\tau$  on the set of all oriented $k$-simplices of $X$ such that
  $ c_{(v_0, \cdots, v_k)}\in {\mathbb C}$;

\noindent (2) $ [v_0, \cdots, v_k]$ is identified with $- [ v_0', \cdots, v_k']$  in the above sum
if $(v_0, \cdots, v_k)$ and $(v'_0, \cdots, v'_k)$ are opposite orientations of the same simplex;

\noindent (3) for any compact subset $K$ of $X$, there are at most finitely many ordered simplices $(v_0, \cdots, v_k) $ intersecting $K$ such that $c_ { (v_0, \cdots, v_k)}\neq 0$.

We remark that in the above definition the summation is taken over all $ (v_0, \cdots, v_k)$ instead of $[v_0, \cdots, v_k]$ for the purpose
to have consistent notations in the definitions of the Connes-Chern characters for lower algebraic $K$-groups of ${\cal S}_1\Gamma$ using simplicial homology groups and  lower algebraic $K$-groups of ${\cal S}_p\Gamma$ using cyclic simplicial homology groups later in this section.

Let $C_k(X)$ be the abelian group of all locally finite $k$-chains on $X$.

Let $$\partial_k: C_k(X) \rightarrow C_{k-1}(X)$$ be the standard simplicial boundary map.
We define the locally finite  simplicial homology group:
$$H_n(X) =Ker~ \partial_n/ Im ~\partial_{n+1}.$$

If $X$ has a proper simplicial action of $\Gamma$,  let $C_k^\Gamma(X)\subset C_k(X)$ be the abelian group consisting of all $\Gamma$-invariant locally finite $k$-chains on $X$.

Let $$\partial^{\Gamma}_k: C^{\Gamma}_k(X) \rightarrow C^{\Gamma}_{k-1}(X)$$ be the restriction of the standard simplicial boundary map.
We define the locally finite $\Gamma$-equivariant  simplicial homology group:
$$H^{\Gamma}_n(X) =Ker~ \partial^{\Gamma}_n/ Im ~\partial^{\Gamma}_{n+1}.$$

Without loss of generality, we can assume that $\Gamma$ is a countable group (this is because every group is an inductive limit of countable groups). We endow $\Gamma$ with a proper left invariant length metric (here properness simply means that every ball with finite radius has finitely many elements). We remark that such a proper  length metric always exists for any countable group.
For each $d\geq 0$, the Rips complex $P_d (\Gamma)$ is the simplicial complex with $\Gamma$  as its vertex set and where  a finite subset $\{\gamma_0, \cdots, \gamma_n\}$  of $\Gamma$ forms a simplex iff $d(\gamma_i, \gamma_j)\leq d$ for all $0 \leq i,j\leq n.$ It is not difficult to show that $\cup_{d\geq 1} P_d (\Gamma)$ is a model for $E_{\cal FIN} (\Gamma)$. When $\Gamma$ is torsion free, $\cup_{d\geq 1} P_d (\Gamma)$ is a universal space for free  and proper action of $\Gamma$. In this case,  $\lim_{d\rightarrow \infty} H_n^{\Gamma} (P_d(\Gamma))$ is the group homology of $\Gamma$ defined using a standard resolution.

To motivate the general construction of the Connes-Chern character, we shall first consider the special case when $\Gamma$ is torsion free.

Let ${\cal A} (\Gamma, {\cal S})$ be the algebra of all kernels
$$k: \Gamma\times \Gamma \rightarrow {\cal S}$$
 such that

\noindent (1) for each $k$, there exists $r\geq 0$ such that $k(x,y)=0$ if $d(x,y) > r$
(the smallest such $r$ is called the propagation of $k$);

\noindent (2) $k$ is $\Gamma$-invariant, i.e. $k(gx, gy)=k(x, y)$ for all $g\in \Gamma$ and $(x,y)\in \Gamma\times \Gamma;$

\noindent (3) The product in ${\cal A} (\Gamma, {\cal S})$  is defined by:
$$(k_1k_2)(x, y)=\sum_{z\in \Gamma} k_1(x, z)k_2 (z, y).$$

We identify ${\cal S}\Gamma$ with ${\cal A} (\Gamma, {\cal S})$ by the isomorphism:
$$\sum_{g\in\Gamma} s_g g\rightarrow k(x, y)=s_{x^{-1}y},$$
where $s_g\in {\cal S}$ for each $g\in \Gamma.$ For each $p\geq 1$, we can naturally identify ${\cal S}_p\Gamma$ with ${\cal A} (\Gamma, {\cal S}_p)$, where ${\cal A} (\Gamma, {\cal S}_p)$ is defined by replacing ${\cal S}$ with ${\cal S}_p$ in the above definition of ${\cal A} (\Gamma, {\cal S})$.

For each non-negative  even integer $n$, we shall first define the Connes-Chern character $c_n$ for a countable torsion free group $\Gamma$
$$c_n: K_0( {\cal S}_1\Gamma) \rightarrow \lim_{d \rightarrow \infty}( \oplus _{k~even, ~k\leq n}~ H_k^{ \Gamma} (P_d (\Gamma)))$$
by:
$$[\tilde{q}]-[q_0]\rightarrow  $$ $$ \sum_{k~ even, ~k\leq n}\sum_{(x_0,\cdots, x_k)}  tr(q(x_0, x_{1})q(x_1, x_{2})\cdots q(x_{k}, x_{0}))[x_{0},\cdots, x_{k}],$$
where
$P_d( \Gamma)$ is the Rips complex and $\tilde{q}$ is an idempotent  in $M_m(({\cal S}_1 \Gamma)^+)$, $\tilde{q}=q+q_0 $ for some $q\in M_m({\cal S}_1\Gamma)$ and idempotent $q_0\in M_m ({\mathbb C})$,
 and the summation is taken over all orderings $(x_0,\cdots, x_k)$  of all $n$-simplices $\{x_0,\cdots, x_k\}$ of $P_d(\Gamma)$ for some $d$ large enough
such that $q(x,y)=0$ if $d(x,y)> d/(n+1)$.

\begin{prop} Let $\Gamma$ be a countable  torsion free group. For each non-negative even integer $n$,
the Connes-Chern character $c_n$ is a well defined \linebreak homomorphism from $ K_0( {\cal S}_1\Gamma)$ to
$\lim_{d \rightarrow \infty}( \oplus _{k~even,~k\leq n}~ H_k^{ \Gamma} (P_d (\Gamma)))$.
\end{prop}
\proof
We first observe that $c( [\tilde{q}]-[q_0])$ is $\Gamma$-invariant by using the $\Gamma$-invariance of $q$.

For each even $k\leq n$, we shall prove that
$$\partial^{\Gamma}_k ( \sum_{(x_0,\cdots, x_k)} tr(q(x_0, x_1)q(x_1, x_2)\cdots q(x_k, x_0))
[x_0,\cdots, x_k] )=0.$$ This implies that $c( [\tilde{q}]-[q_0])$ is a cycle.

 We leave  to the reader the proof that the homology class of
 $$ \sum_{k~ even,~k\leq n}\sum_{(x_0,\cdots, x_k)}  tr(q(x_0, x_{1})q(x_1, x_{2})\cdots q(x_{k}, x_{0}))[x_{0},\cdots, x_{k}]$$
depends only on the K-theory class $[\tilde{q}]-[q_0]$.

By the assumption that $\tilde{q}$ and $q_0$ are idempotents, we have
$$q^2=q-q_0q-qq_0.$$

It follows that
\begin{eqnarray*}
& &\partial_k ( \sum_{(x_0,\cdots, x_k)} tr(q(x_0, x_1)q(x_1, x_2)\cdots q(x_k, x_0))
[x_0,\cdots, x_k] )= \\
 &&      \sum_{i=1}^k    (-1)^i \sum_{(x_0,\cdots,\hat{x_i},\cdots, x_k)}
 tr(q(x_0, x_1)q(x_1, x_2)\cdots q^2(x_{i-1}, x_{i+1}) \cdots q(x_k, x_0))
[x_0,\cdots, \hat{x_i}, \cdots, x_k]                                  \\
& & + \sum_{(x_1,\cdots, x_k)} tr(q(x_1, x_2)\cdots q^2(x_k, x_1))
[x_1,\cdots, x_k]= \\
\end{eqnarray*}

\begin{eqnarray*}
&&  \sum_{i=1}^k   (-1)^i  \sum_{(x_0,\cdots,\hat{x_i},\cdots, x_k)}
 tr(q(x_0, x_1)q(x_1, x_2)\cdots q(x_{i-1}, x_{i+1}) \cdots q(x_k, x_0))
[x_0,\cdots, \hat{x_i}, \cdots, x_k]                                  \\
& & + \sum_{(x_1,\cdots, x_k)} tr(q(x_1, x_2)\cdots q(x_k, x_1))
[x_1,\cdots, x_k]-\\
& & (\sum_{i=1}^k   (-1)^i  \sum_{(x_0,\cdots,\hat{x_i},\cdots, x_k)}
 tr(q(x_0, x_1)q(x_1, x_2)\cdots q_0q(x_{i-1}, x_{i+1}) \cdots q(x_k, x_0))
[x_0,\cdots, \hat{x_i}, \cdots, x_k]                                  \\
& & + \sum_{(x_1,\cdots, x_k)} tr(q(x_1, x_2)\cdots q_0q(x_k, x_1))
[x_1,\cdots, x_k] +\\
& &   \sum_{i=1}^k   (-1)^i \sum_{(x_0,\cdots,\hat{x_i},\cdots, x_k)}
 tr(q(x_0, x_1)q(x_1, x_2)\cdots q(x_{i-1}, x_{i+1})q_0 \cdots q(x_k, x_0))
[x_0,\cdots, \hat{x_i}, \cdots, x_k]                                  \\
& & + \sum_{(x_1,\cdots, x_k)} tr(q(x_1, x_2)\cdots q(x_k, x_1)q_0)
[x_1,\cdots, x_k] ).
\end{eqnarray*}
Using the trace property and the definition of oriented simplices and the assumption that $k$ is even, we have
\begin{eqnarray*}
& &\sum_{(x_0,\cdots,\hat{x_i},\cdots, x_k)}
 tr(q(x_0, x_1)q(x_1, x_2)\cdots q(x_{i-1}, x_{i+1}) \cdots q(x_k, x_0))
[x_0,\cdots, \hat{x_i}, \cdots, x_k]=0 , \\
&  & \sum_{(x_1,\cdots, x_k)} tr(q(x_1, x_2)\cdots q(x_k, x_1))
[x_1,\cdots, x_k] =0,\\
& & \sum_{(x_0,\cdots,\hat{x_i},\cdots, x_k)}
 ( tr(q(x_0, x_1)q(x_1, x_2)\cdots q_0q(x_{i-1}, x_{i+1}) \cdots q(x_k, x_0))
[x_0,\cdots, \hat{x_i}, \cdots, x_k] + \\
& & tr(q(x_0, x_1)q(x_1, x_2)\cdots q(x_{i-1}, x_{i+1})q_0 \cdots q(x_k, x_0))
[x_0,\cdots, \hat{x_i}, \cdots, x_k]  )=0,\\
& &\sum_{(x_1,\cdots, x_k)}( tr(q(x_1, x_2)\cdots q_0q(x_k, x_1))
[x_1,\cdots, x_k] + \\
& & tr(q(x_1, x_2)\cdots q(x_k, x_1)q_0)
[x_1,\cdots, x_k] )=0.
\end{eqnarray*} \qed

By the definition of lower algebraic K-theory using the group algebra over the free abelian group $ \mathbb{Z}^n$, we can similarly define the Connes-Chern character:
$$c_n: K_i( {\cal S}_1\Gamma) \rightarrow \lim_{d \rightarrow \infty} (\oplus _{k+i~even, ~k\leq n} ~H_k^{ \Gamma} (P_d (\Gamma)))$$
for each non-negative  integer $n$ and $i<0$.

Next we extend the construction of the  Connes-Chern character $c_n$ to the ${\cal S}_p$ case for each $p\geq 1$ and each non-negative integer
$n$ when $\Gamma$ is torsion free:
$$c_n: K_0( {\cal S}_p\Gamma) \rightarrow \lim_{d \rightarrow \infty}( \oplus _{k~even,~k\leq n}~ H_k^{\Gamma} (P_d (\Gamma))).$$

We need to introduce an equivariant  cyclic simplicial homology group to define the Connes-Chern character.

Let $X$ be a simplicial complex. An ordered $k$-simplex $(v_0, \cdots, v_k)$ is defined to be  an  ordered  finite sequence of vertices in a simplex of $X$,
where $v_i$ is allowed to be equal to $v_j$ for some distinct pair of $i$ and $j$.

Recall that the following permutation  is called a cyclic permutation
$$(v_0, ..., v_k)\rightarrow (v_k,  v_0, \cdots, v_{k-1}).$$
 We define two ordered simplices $(v_0, \cdots, v_k)$  and $(v'_0, \cdots, v'_k)$  to be equivalent if one ordered  simplex can be obtained from the other ordered simplex by any number of  cyclic permutations when $k$ is even and  by an even number of  cyclic permutations when $k$ is odd. Each of the equivalence classes is called a cyclically oriented simplex.
 If $(v_0, ..., v_k)$ is an ordered simplex of $X$, we use the symbol
$[v_0, \cdots, v_k]_\lambda$ to denote the corresponding cyclically oriented simplex.

A locally finite cyclic $k$-chain on $X$ is  a formal sum
$$\sum_{(v_0, \cdots, v_k)} c_{(v_0, \cdots, v_k)} [v_0, \cdots, v_k]_{\lambda},$$ where
%${\mathbb R}$-valued  function $\tau$  on the set of all oriented $k$-simplices of $X$ such that

\noindent (1) the summation is taken over all  ordered simplices $ (v_0, \cdots, v_k)$ of   $X$ and  $ c_{(v_0, \cdots, v_k)}\in {\mathbb C}$;

\noindent (2) $ [v_0, \cdots, v_k]_{\lambda}$ is identified with $(-1)^k [v_k, v_0, \cdots, v_{k-1}]_{\lambda}$ in the above sum;

\noindent (3)
for any compact subset $K$ of $X$, there are at most finitely many ordered  simplices $(v_0, \cdots, v_k) $ intersecting $K$ such that $c_{(v_0, \cdots, v_k)}\neq 0$.

Let $C^{\lambda}_k(X)$ be the abelian group of all locally finite cyclic $k$-chains on $X$.
Let $$\partial^{\lambda}_k: C^{\lambda}_k(X) \rightarrow C^{\lambda}_{k-1}(X)$$
 be the standard boundary map.

We define the cyclic simplicial homology group:
$$H_n^{\lambda}(X) =Ker~ \partial_n^{\lambda}/ Im ~\partial_{n+1}^{\lambda}.$$

If $X$ has a simplicial  proper action of a group $\Gamma$, we can  define $C^{\lambda, \Gamma}_k(X)\subseteq C^{\lambda}_k(X)$ to be the subspace of all $\Gamma$-invariant locally finite cyclic $k$-chains on $X$.

Let $$\partial^{\lambda,\Gamma}_k: C^{\lambda,\Gamma}_k(X) \rightarrow C^{\lambda,\Gamma}_{k-1}(X)$$ be the restriction of the standard boundary map.

We define the $\Gamma$-equivariant cyclic simplicial homology group
$H_n^{\lambda, \Gamma}(X)$ by:
$$H_n^{\lambda,\Gamma}(X) =Ker~ \partial_n^{\lambda,\Gamma}/ Im ~\partial_{n+1}^{\lambda,\Gamma}.$$

The following result computes the $\Gamma$-equivariant cyclic simplicial homology group in terms of
the $\Gamma$-equivariant  simplicial homology groups.

\begin{prop} Let $\Gamma$ be a group.  Let $X$ be a simplicial complex with a proper simplicial action of $\Gamma$.
We have $$ H_n ^{\lambda, \Gamma}(X)~\cong~~ (\oplus_{k\leq n,~ k=n~mod~2}~ H_k^{\Gamma}(X)).$$
\end{prop}
\proof
Given an ordering $(v_0, \cdots, v_k)$ of $k+1$ number of  vertices in $X$, we use the same notation  $(v_0, \cdots, v_k)$ to  denote the corresponding ordered $k$-simplex.  Let $C^{ord, \Gamma}_{ k}(X)$ be the abelian group  of all $\Gamma$-invariant
locally finite ordered $k$-chains
$$\sum_{(v_0, \cdots, v_k)} c_{(v_0, \cdots, v_k)}(v_0, \cdots, v_k),$$
where the sum is taken over all ordered $k$-simplices of $X$ and $c_{(v_0, \cdots, v_k)}\in \mathbb{C}.$
We remark that, in the above definition,  a pair of vertices  in $(v_0, \cdots, v_k)$ are allowed to be the same.

Let $C^{ord, \Gamma}_{ k,0}(X)$ be the abelian subgroup of $C^{ord, \Gamma}_{ k}(X)$  consisting of all $\Gamma$-invariant
locally finite ordered $k$-chains with the following special form:
$$\sum_{v} c_{v}\overbrace{(v, \cdots, v)}^\text{k+1},$$
where the sum is taken over all vertices of $X$ and $c_v\in \mathbb{C}.$

We define an abelian group  $C^{ord, \Gamma}_{ k, red}(X)$ by:
$$
 C^{ord, \Gamma}_{ k, red}(X)  := \left\{\begin{array}{ll} C^{ord, \Gamma}_{ k}(X)  &  \text{if $k$ is even,} \\ C^{ord, \Gamma}_{ k}(X)/C^{ord, \Gamma}_{ k,0}(X) & \text{if $k$ is odd,}\end{array}\right.
$$ where $C^{ord, \Gamma}_{ k}(X)/C^{ord, \Gamma}_{ k,0}(X)$ is the quotient group of $C^{ord, \Gamma}_{ k}(X)$ over $C^{ord, \Gamma}_{ k,0}(X)$.

The standard  boundary map on $C^{ord, \Gamma}_{ k}(X)$ induces a boundary map: $$\partial^{ord, \Gamma}_{k, red}:
  C^{ord, \Gamma}_{ k, red}(X) \longrightarrow C^{ord, \Gamma}_{ k-1, red}(X)  .$$

We define a new homology group:
$$H_{n, red}^{ord, \Gamma}(X) =Ker~ \partial^{ord,\Gamma}_{n, red}/ Im ~\partial^{ord,\Gamma}_{n+1, red}.$$

Let $\chi_k$ be the natural chain map from $C^{ord, \Gamma}_{ k,red }(X)$ to $C^{\lambda, \Gamma}_{ k}(X)$ defined by:
$$[\sum_{(v_0, \cdots, v_k)} c_{(v_0, \cdots, v_k)}(v_0, \cdots, v_k)] \longrightarrow  \sum_{(v_0, \cdots, v_k)} c_{(v_0, \cdots, v_k)}[v_0, \cdots, v_k]_{\lambda},$$
for every  $[\sum_{(v_0, \cdots, v_k)} c_{(v_0, \cdots, v_k)}(v_0, \cdots, v_k)] \in C^{ord, \Gamma}_{ k,red }(X)$. This map is well defined because
$\overbrace{[v, \cdots, v]}^\text{k+1}$ $_\lambda=0$ when $k$ is odd.

By a standard Mayer-Vietoris  and five lemma argument,
it is not difficult to prove that $\chi$ induces an isomorphism $\chi_\ast$ from $H_{n, red}^{ord, \Gamma}(X) $ to $H_{n}^{ \lambda, \Gamma}(X)$.
This is because both homology theories satisfies the Mayer-Vietoris sequences and $\chi_\ast$ is an isomorphism when $X=\Gamma/F$ as $\Gamma$-spaces for some finite subgroup $F$ of $\Gamma$.

We define a natural chain map $$\phi_{k,k}: C^{ord,\Gamma}_{ k, red}(X)\rightarrow C^{\Gamma}_{k}(X)$$  by:
$$[\sum_{(v_0, \cdots, v_k)} c_{(v_0, \cdots, v_k)}(v_0, \cdots, v_k)] \longrightarrow  \sum_{(v_0, \cdots, v_k)} c_{(v_0, \cdots, v_k)}[v_0, \cdots, v_k]$$
for every  $[\sum_{(v_0, \cdots, v_k)} c_{(v_0, \cdots, v_k)}(v_0, \cdots, v_k)] \in C^{ord, \Gamma}_{ k,red }(X)$. This map is well defined because
$\overbrace{[v, \cdots, v]}^\text{k+1}=0$ when $k$ is odd (more generally when $k\geq 1$).

 For each  ordered $k$-simplex $ (v_0, \cdots, v_k)$ of $X$, we let
$$
 \varphi_{k,k-2}( (v_0, \cdots, v_k))  :=$$
  $$\left\{\begin{array}{ll} (-1)^{j-i+1}[v_0, \cdots, \hat{v_i}, \cdots, \hat{v_j}, \cdots, v_k]  &  \text{if $(i,j)$ is the smallest pair such that $i<j$, $v_i=v_j$ ,}
  \\ $~~~~~~~~~~~~~~0$ & \text{if there exists no pair $i<j$ such that $v_i=v_j$,}\end{array}\right.
$$  
  where the smallest $(i,j)$ is taken with respect to the dictionary order of $\{(m,l): 0\leq m<l\leq k\}$ given by:
 $(m, l)< (m', l')$ iff either (1) $ m<m'$,  or (2) $m=m'$ and $l<l'$.

We define a linear map $$\phi_{k,k-2}: C^{ord,\Gamma}_{ k, red}(X)\rightarrow C^{\Gamma}_{k-2}(X)$$  by:
   $$[\sum_{(v_0, \cdots, v_k)} c_{(v_0, \cdots, v_k)}(v_0, \cdots, v_k)] \longrightarrow  
   \sum_{(v_0, \cdots, v_k)} c_{(v_0, \cdots, v_k)}\varphi_{k,k-2}((v_0, \cdots, v_k))$$
for every  $[\sum_{(v_0, \cdots, v_k)} c_{(v_0, \cdots, v_k)}(v_0, \cdots, v_k)] \in C^{ord, \Gamma}_{ k,red }(X)$.                  
 
Note that $\phi_{k,k-2}$ is well defined. Elementary computations show that $\phi_{k,k-2}$ is a chain map.

Similarly we can construct a chain map
$$\phi_{k, l}: C^{ord,\Gamma}_{ k, red}(X)\rightarrow C^{\Gamma}_{ l}(X)$$ if $0\leq l\leq k$ and $k-l$ is even.

Using the above chain maps, we  construct a chain map:
$$\psi_n=\oplus_{ l\leq n,~ n-l ~even}~ \phi_{n, l} : C^{ord,\Gamma}_{n, red}(X) \rightarrow (\oplus_{l \leq n,~ n-l~even}~C^{\Gamma}_{l}(X)).$$

The   chain map $\psi_n$  induces an isomorphism $\psi_\ast$ on the homology groups  since both homology theories satisfy the Mayer-Vietoris sequence
and the chain map induces an isomorphism at the homology level if $X=\Gamma/F$ as $\Gamma$-spaces for some finite subgroup $F$ of $\Gamma$.

Finally Proposition 4.2 follows from the facts that $\chi_\ast$ and $\psi_\ast$ are isomorphisms.
\qed

%Define a map $$ S: ~~~~\overbrace{ {\cal S}_p\otimes \cdots \otimes {\cal S}_p}^{n+3} \longrightarrow \overbrace{ {\cal S}_p\otimes \cdots \otimes {\cal S}_p}^{n+1}$$ by

%\begin{eqnarray*}
%S(a_0\otimes  a_1\otimes \cdots\otimes a_{n+2})  & = &  (a_0a_1 a_2)\otimes a_3\otimes\cdots\otimes a_{n+2}+ \\
% & &  ( a_0\otimes (a_1a_2 a_3)\otimes \cdots \otimes a_{n+2} -\\
% & & (a_0a_1)\otimes (a_2a_3)\otimes a_4 \otimes \cdots \otimes a_{n+2}) +\\
% & & ( a_0 \otimes a_1\otimes \cdots \otimes (a_n a_{n+1}a_{n+2})-\\
% & & a_0 \otimes \cdots \otimes (a_{n-1}a_n) \otimes (a_{n+1}a_{n+2})+\\
% & & a_0\otimes \cdots \otimes (a_{n-2}a_{n-1})\otimes a_n \otimes (a_{n+1}a_{n+2})-\\
% & &\cdots + \\
% & & (-1)^n (a_0a_1a_2) \otimes a_3 \otimes \cdots \otimes a_n \otimes (a_{n+1}a_{n+2})).
% \end{eqnarray*}

For any positive even integer $n\geq p$, we  now define the Connes-Chern character for a countable torsion free group $\Gamma$
$$c_n: K_0( {\cal S}_p\Gamma) \rightarrow H_n^{\lambda, \Gamma} (P_d (\Gamma)),$$
where $\Gamma$ is endowed with a proper length metric.

Let $\tilde{q}$ be an idempotent  in $M_m(({\cal S}_p \Gamma)^+)$ and $\tilde{q}=q+q_0$ for some $q\in M_m({\cal S}_p\Gamma)$ and    idempotent $q_0\in M_m ({\mathbb C})$. We  identify $q$ with an element in ${\cal A} (\Gamma, {\cal S}_p)$
(note that ${\cal A} (\Gamma, M_m( {\cal S}_p))$ is isomorphic to  ${\cal A} (\Gamma, {\cal S}_p)$).
 Let $d$ be greater than or equal to $n+1$ times the propagation of $q$, i.e. $q(x,y)=0$ if $d(x,y)>d/(n+1)$.

For each positive even integer $n\geq p$, the Connes-Chern character $c_n$ of
$[\tilde{q}] -[q_0]$ is defined to be homology class of
$$\sum_{(x_0,\cdots, x_n) } tr(q(x_0, x_1) \cdots q(x_n, x_0))
[x_0,\cdots, x_n]_{\lambda}\in H_n^{\lambda, \Gamma} (P_d (\Gamma)),$$
where the summation is taken over all ordered $n$-simplices $(x_0,\cdots, x_n)$ of $P_d(\Gamma)$.

We remark that the choice of $n$ guarantees that the trace  in the above definition of the Connes-Chern character is finite.

By Proposition 4.2, the above Connes-Chern character induces a Connes-Chern character:

$$c_n:  K_0( {\cal S}_p\Gamma)\rightarrow
\lim_{d \rightarrow \infty}( \oplus _{k~even,~k\leq n}~ H_k^{ \Gamma} (P_d (\Gamma)))$$
for any non-negative integer $n\geq p$.

For an arbitrary non-negative even integer $n$, let $n'$ be a positive even integer satisfying $n'\geq max\{n, p\}.$
Let $\pi_{n', n}$ be the natural projection from \linebreak $\lim_{d \rightarrow \infty}( \oplus _{k~even, ~k\leq n'} ~H_{k}^{ \Gamma} ( P_d (\Gamma)))$
to $\lim_{d \rightarrow \infty}( \oplus _{k~even, ~k\leq n} ~H_{k}^{ \Gamma} ( P_d (\Gamma)))$. We define the Connes-Chern character
$c_n$ from $ K_0( {\cal S}_p\Gamma)$ to \linebreak
$\lim_{d \rightarrow \infty}( \oplus _{k~even, ~k\leq n} ~ H_{k}^{ \Gamma} ( P_d (\Gamma)))$ to be $\pi_{n', n} \circ c_{n'}$. It is not difficult to verify that the definition of $c_n$ is independent of the choice of $n'$.

\begin{prop} Let $\Gamma$ be a countable  torsion free group. For any non-negative even integer $n$,
the Connes-Chern character $c_n$ is a well defined homomorphism from $ K_0( {\cal S}_p\Gamma)$ to
$\lim_{d \rightarrow \infty}( \oplus _{k~even,~k\leq n}~ H_k^{ \Gamma} (P_d (\Gamma)))$.
\end{prop}

The proof of the above proposition is similar to the proof of Proposition 4.1 and is therefore omitted.
Note that when $p=1$, the above definition of the Connes-Chern character coincides with the prior definition of the Connes-Chern character.

Next we shall construct the Connes-Chern character for a general group  $\Gamma$.

Let $\Gamma_{fin}$ be the set of all elements with finite order in $\Gamma$.  The group $\Gamma$ acts on $\Gamma_{fin}$ by conjugations:
$$\gamma\cdot x =\gamma x \gamma^{-1}$$ for all $\gamma \in \Gamma$ and $x\in \Gamma_{fin}.$

Let   $X$ be a simplicial complex with a proper simplicial action of $\Gamma$. Equip the vertex set $V(X)$ of $X$ with a $\Gamma$-invariant proper pseudo metric $d_V$.
Let $\Gamma$ act on $\Gamma_{fin} \times X$ diagonally.

Let $r\geq 0$. For each $g\in \Gamma_{fin}$, we define
$X_{g, r} $ to be the simplicial subcomplex  of $X$ consisting  all simplices $\{v_0, \cdots, v_p\}$  satisfying $ d_V (v_i, gv_i)\leq r$ for all $0\leq i\leq p.$

For each ordered simplex $(v_0, ..., v_k)$ of $X_{g, r}$, we define the following transformation to be a $g$-cyclic permutation
$$(v_0, ..., v_k)\rightarrow (gv_k,  v_0, \cdots, v_{k-1}).$$
We define two ordered simplices $(v_0, \cdots, v_k)$  and $(v'_0, \cdots, v'_k)$ of $X_{g, r}$ to be $g$-equivalent if one ordered  simplex can be obtained from the other ordered simplex by any number of  $g$-cyclic permutations  of ordered simplices in $X_{g, r}$  when $k$ is even and  by an even number of  $g$-cyclic permutations when $k$ is odd. Each of the equivalence classes is called a $g$-cyclically oriented simplex.
 If $(v_0, ..., v_k)$ is an ordered simplex of $X_{g, r}$, we use the symbol
$[v_0, \cdots, v_k]_{\lambda,g}$ to denote the corresponding $g$-cyclically oriented simplex.

We define ${\mathbb C}^{\lambda}_{k,r}(X) $ to be the abelian group  of all locally finite  $k$-chains:
$$ \sum_{g\in \Gamma_{fin} }  (g, \sum_{(v_0, \cdots, v_k)} c_{(v_0, \cdots, v_k),g} [v_0, \cdots, v_k]_{\lambda, g}),$$ where

\noindent (1) the second summation is taken over all  ordered simplices $ (v_0, \cdots, v_k)$ of   $X_{g, r}$ and  $ c_{(v_0, \cdots, v_k), g}\in {\mathbb C}$;

\noindent (2) $ [v_0, \cdots, v_k]_{\lambda, g}$ is identified with $(-1)^k [gv_k, v_0, \cdots, v_{k-1}]_{\lambda,g}$ in the above sum;

\noindent (3)
for each $g\in \Gamma_{fin}$ and any compact subset $K$ of $X$, there are at most finitely many ordered  simplices $(v_0, \cdots, v_k) $ intersecting $K$ such that $c_{(v_0, \cdots, v_k),g}\neq 0$.

The diagonal action of $\Gamma$ on $\Gamma_{fin}\times X$ induces
 a natural $\Gamma$-action on ${\mathbb C}^{\lambda}_{k,r}(X) $. Let ${\mathbb C}_{k,r}^{\lambda, \Gamma}(X) \subseteq {\mathbb C}_{k,r}^{\lambda}(X)  $ be the abelian group consisting of all  $\Gamma$-invariant $k$-chains in ${\mathbb C}^{\lambda}_{k,r}(X) $.

We have a natural boundary map:
 $$\partial_{k,r}^{\lambda, \Gamma}: {\mathbb C}_{k,r}^{\lambda, \Gamma}(X) \longrightarrow {\mathbb C}_{k-1,r}^{\lambda, \Gamma}(X) .$$
We define the following equivariant homology theory by:

$$ {\mathbb H}^{\lambda, \Gamma}_{n,r} (X)= Ker~\partial_{n,r}^{\lambda, \Gamma}/Im~\partial_{n+1,r}^{\lambda, \Gamma}.$$

When $\Gamma$ is torsion free, $\Gamma_{fin}$ consists of the identity element and  we have
 $$ {\mathbb H}^{\lambda, \Gamma}_{n,r} (X)= H^{\lambda, \Gamma}_n(X).$$

For each $r\geq 0$, let  $\hat{X}_r$ be the simplical subspace of $\Gamma_{fin}\times X$ defined by:
$$\hat{X}_r=\{ (g, x)\in \Gamma_{fin}\times X: x\in X_{g,r}\}.$$
 The diagonal action of $\Gamma$ on $\Gamma_{fin}\times X$ induces
 a natural $\Gamma$-action on $\hat{X}_r$.

We define
$$ {\mathbb H}_{n,r} ^{\Gamma}(X)= H_n^{\Gamma}( \hat{X}_r ).$$

The following result computes our new equivariant homology theory of the Rips complex in terms of the (locally finite) equivariant simplicial homology theory.

\begin{prop} Let $\Gamma$ be a countable group with a proper length metric. We have $$ \lim_{d\rightarrow\infty}\lim_{r\rightarrow\infty}{\mathbb H}_{n,r} ^{\lambda, \Gamma}(P_d(\Gamma))~\cong~~ \lim_{d\rightarrow\infty}\lim_{r\rightarrow \infty}(\oplus_{k\leq n,~ k=n~mod~2}~ {\mathbb H}_{n,r} ^{\Gamma}(P_d(\Gamma))).$$
\end{prop}
\proof
Let $X$ be a simplicial complex with a proper and cocompact action of $\Gamma$. We define an equivalence relation $\sim$
 on the chain group ${\mathbb C}^{\lambda, \Gamma}_{k,r}(X ) $ as follows.
Two chains $z$ and $z'$  in ${\mathbb C}^{\lambda, \Gamma}_{k,r}(X ) $ are said to be equivalent if
$$ z= \sum_{g\in \Gamma_{fin} }  (g, \sum_{(v_0, \cdots, v_k)} c_{(v_0, \cdots, v_k),g} [v_0, \cdots, v_k]_{\lambda, g}),$$
$$ z' =
\sum_{g\in \Gamma_{fin} }  (g, \sum_{(v_0, \cdots, v_k)} c_{(v_0, \cdots, v_k),g} [v'_0, \cdots, v'_k]_{\lambda, g}),$$
and for each $0\leq i \leq k$ there exists an integer $j$ such that  $v'_i=g^jv_i$.

Let $C_{k, r}^{\lambda, \Gamma}(X)$ be the chain group ${\mathbb C}^{\lambda, \Gamma}_{k,r}(X )/\sim$. We define $  \tilde{{\mathbb H}}_{n, r}^{\lambda,\Gamma} (X)$  to be the n-th homology group of
$C_{k, r}^{\lambda, \Gamma}(X)$.

 The quotient chain map $\phi$ from ${\mathbb C}^{\lambda, \Gamma}_{k,r}(X ) $ to $C_{k, r}^{\lambda, \Gamma}(X)$ induces a homomorphism
$$\phi_\ast: {\mathbb H}_{n,r} ^{\lambda, \Gamma}(X)\rightarrow \tilde{ {\mathbb H}}_{n, r}^{\lambda, \Gamma} (X).$$

 We observe that the cocompactness of the $\Gamma$ action on $X$ implies that, for each $r\geq 0$, there exists $N>0$ such that if $g\in \Gamma_{fin}$ and $X_{g, r}$ is nonempty, then
the order of the group element $g$  is  bounded by $N$.
As a consequence, for any $d\geq 0$ and $r\geq 0$, there exist $d' \geq d$ and $r'\geq r$ such that, for any $g\in\Gamma_{fin}$ and any simplex in $(P_d(\Gamma))_{g, r}$ with vertices $\{ v_0, \cdots, v_k\}$, the simplex with vertices $\{ g^{i_0}v_0, \cdots, g^{i_k}v_k: 1\leq i_j\leq N, 0\leq j\leq k\}$ is a simplex in $(P_{d'}(\Gamma))_{g, r'}$.
This implies that $\phi$ is a chain homotopy equivalence
from the chain complex $\lim_{d\rightarrow\infty} \lim_{r\rightarrow\infty} {\mathbb C}^{\lambda, \Gamma}_{k,r}(P_d(\Gamma) ) $
to the chain complex $\lim_{d\rightarrow\infty} \lim_{r\rightarrow\infty} C^{\lambda, \Gamma}_{k,r}(P_d(\Gamma) ) $ with a homotopy inverse chain map $\psi$ from $\lim_{d\rightarrow\infty} \lim_{r\rightarrow\infty} C^{\lambda, \Gamma}_{k,r}(P_d(\Gamma) ) $ to $\lim_{d\rightarrow\infty} \lim_{r\rightarrow\infty} {\mathbb C}^{\lambda, \Gamma}_{k,r}(P_d(\Gamma) ) $  defined by
\begin{eqnarray*}
& & \psi([\sum_{g\in \Gamma_{fin} }  (g, \sum_{(v_0, \cdots, v_k)} c_{(v_0, \cdots, v_k),g} [v_0, \cdots, v_k]_{\lambda, g})]) = \\
 &&      \sum_{g\in \Gamma_{fin} }  (g, \sum_{(v_0, \cdots, v_k)}\frac{1}{n_g^{k+1}}\sum_{ i_0, \cdots, i_k=1  }^{n_g} c_{(v_0, \cdots, v_k),g} [g^{i_0}v_0, \cdots, g^{i_k}v_k]_{\lambda, g})
\end{eqnarray*}
for each  $[\sum_{g\in \Gamma_{fin} }  (g, \sum_{(v_0, \cdots, v_k)} c_{(v_0, \cdots, v_k),g} [v_0, \cdots, v_k]_{\lambda, g})]$ in \linebreak
$\lim_{d\rightarrow\infty} \lim_{r\rightarrow\infty} C^{\lambda, \Gamma}_{k,r}(P_d(\Gamma) ) $,
where $n_{g}$ is the order of the group element $g$.
It follows  that
the homomorphism $\phi_\ast$ is an isomorphism from \linebreak $\lim_{d\rightarrow\infty} \lim_{r\rightarrow\infty}{\mathbb H}_{n,r} ^{\lambda, \Gamma}(P_d(\Gamma))$
to $\lim_{d\rightarrow\infty} \lim_{r\rightarrow\infty}  \tilde{ {\mathbb H}}_{n,r} ^{\lambda, \Gamma}(P_d(\Gamma)).$

Two vertices $v$ and $v'$ of $X_{g,r}$ are defined to be equivalent if $v=g^j v'$ for some $j$.
We denote the equivalence class of $v$ by $[v]$.

We define $\tilde{X}_{g,r}$ to be the simplicial complex consisting of simplices $\{[v_0], \cdots, [v_k]\}$  for all simplices
$\{v_0, \cdots, v_k\} $ in $X_{g,r}$.

Let $$ \tilde{X}_r= \{(g, x): g\in \Gamma_{fin}, x\in \tilde{X}_{g, r}\} .$$

By an argument similar to the proof of Proposition 4.2, we have the following isomorphism:
$$ \lim_{d\rightarrow\infty}\lim_{r\rightarrow\infty}  \tilde{ {\mathbb H}}_{n,r} ^{\lambda, \Gamma}(P_d(\Gamma)) \cong
\lim_{d\rightarrow\infty}\lim_{r\rightarrow \infty}(\oplus_{k\leq n,~ k=n~mod~2}~  H_{n} ^{\Gamma}((\widetilde{P_d(\Gamma)})_r)).$$

Finally we observe that the natural homomorphism
$$\lim_{d\rightarrow\infty}\lim_{r\rightarrow \infty}(\oplus_{k\leq n,~ k=n~mod~2}~ {\mathbb H}_{n,r} ^{\Gamma}(P_d(\Gamma)))\rightarrow
\lim_{d\rightarrow\infty}\lim_{r\rightarrow \infty}(\oplus_{k\leq n,~ k=n~mod~2}~  H_{n} ^{\Gamma}((\widetilde{P_d(\Gamma)})_r))$$
is an isomorphism.
\qed

Let $\Gamma$ be a countable group with a proper length metric.  Let $X$ be a simplicial complex with a proper and cocompact action of $\Gamma$.
 Let  $\hat{X}$ be the  subspace of $\Gamma_{fin}\times X$ defined by:
$$\hat{X}=\{ (g, x)\in \Gamma_{fin}\times X: gx=x \}.$$
 The diagonal action of $\Gamma$ on $\Gamma_{fin}\times X$ induces
 a natural $\Gamma$-action on $\hat{X}$. Note that $\hat{X}$ is a simplicial complex with a simplicial action of $\Gamma$.

 We define
 $${\mathbb H}_k^{ \Gamma} (X)= H^{\Gamma}_k (\hat{X}).$$

 We remark that ${\mathbb H}_k^{ \Gamma} (X)$ is the equivariant homology theory of Baum-Connes [BC2].

\begin{prop} Let $\Gamma$ be a countable group with a proper length metric.
We have $$\lim_{d\rightarrow\infty} \lim_{r\rightarrow \infty}{\mathbb H}_{n,r} ^{\Gamma}(P_d(\Gamma))\cong \lim_{d\rightarrow\infty}
{\mathbb H}_{n}^{\Gamma}(P_d(\Gamma)).$$
\end{prop}
\proof
Let $X$ be a simplicial complex with a proper and cocompact action of $\Gamma$.
For each finite subset $F\subset \Gamma_{fin}$, let
$$F'=\{\gamma f\gamma^{-1}: \gamma\in \Gamma, f\in F\}.$$

For each $g\in \Gamma_{fin}$ and $r\geq 0$, we define $\check{X}_{g, r}$ to be the simplicial subcomplex of $X$ consisting of  simplices with vertices
$$\{g^{j_0} v_0, \cdots, g^{j_k}v_k: j_i \in {\mathbb Z} \}$$ for  all simplices $\{v_0, \cdots, v_k\}$ in $X_{g,r}$.

We let $$\hat{X}_{F}=\{ (g, x)\in F'\times X: gx=x\}, ~~~~~\check{X}_{F, r}= \{ (g, x)\in F'\times X: x\in \check{X}_{g, r}\}.$$

We have an inclusion map
$$i: \hat{P_d (\Gamma)}_F \rightarrow \check{(P_d(\Gamma))}_{F, r}.$$
The map $i$ induces a homomorphism
$$ i_\ast: \lim_{d\rightarrow\infty}{\mathbb H}_{n}^{\Gamma}(P_d(\Gamma)) \rightarrow \lim_{d\rightarrow\infty} \lim_{r\rightarrow \infty}{\mathbb H}_{n,r} ^{\Gamma}(P_d(\Gamma)).$$

By the definition of $\check{(P_d(\Gamma))}_{F,r}$, for
 each $d\geq 0$ and $r\geq 0$,  there exists $c>0$ such that, for every point $(g,x)$ in $\check{(P_d(\Gamma))}_{F,r}$, $x$ is within  distance $c$ from a fixed point of $g$.
It follows that, for each $d\geq 0$ and $r\geq 0$, there exist $d'\geq d$ and a continuous map
$$\psi: \check{(P_d(\Gamma))}_{F, r} \rightarrow \hat{(P_{d'} (\Gamma))}_F$$
such that if we write
 $\psi (g, x)= (g, \psi' (x))$, then we have
 $$\sup \{ d( \psi'(x), x): (g,x) \in \hat{(P_{d'} (\Gamma))}_F\}<\infty,$$
 where $d$ is the restriction of the simplicial metric on $P_d(\Gamma)$.

The map $\psi$ induces a homomorphism
 $$ \psi_\ast: \lim_{d\rightarrow\infty} \lim_{r\rightarrow \infty}{\mathbb H}_{n,r} ^{\Gamma}(P_d(\Gamma)) \rightarrow   \lim_{d\rightarrow\infty}{\mathbb H}_{n}^{\Gamma}(P_d(\Gamma)).$$

 Using linear homotopies, it is not difficult to check that $i_\ast$ and $\psi_\ast$  are inverses to each other.
\qed

For each non-negative integer $n$, we are now ready to define the Connes-Chern character $c_n$
 for a general group $\Gamma$:
$$c_n:~~ K_0( {\cal S}\Gamma) \longrightarrow  \lim_{d\rightarrow \infty} (\oplus _{k~even,~k\leq n} ~~ {\mathbb H}_{n}^{ \Gamma} ( P_d(\Gamma))).$$

For each $p\geq 1$ and  even integer $n\geq p$, we shall first define the Connes-Chern character $c_m$:
$$c_n:~~ K_0( {\cal S}_p\Gamma) \longrightarrow   \lim_{d\rightarrow \infty} (\oplus _{k~even, ~k\leq n} ~~ {\mathbb H}_{k}^{ \Gamma} ( P_d(\Gamma)))$$

Let $\tilde{q}$ be an idempotent  in $M_m(({\cal S}_p \Gamma)^+)$ and $\tilde{q}=q+q_0$ for some $q\in M_m({\cal S}_p\Gamma)$ and    idempotent $q_0\in M_m ({\mathbb C})$. Let $d$ be greater than or equal to $n+1$ times  the propagation of $q$, i.e. $q(x,y)=0$ if $d(x,y)>d/(n+1)$.

The Connes-Chern character $c_n$ of
$[\tilde{q}] -[q_0]$ is defined to be homology class of
$$\sum_{g\in \Gamma_{fin}} (g, \sum_{(x_0,\cdots, x_n) } tr(q(x_0, x_1) \cdots q(x_n, g^{-1}x_0))
[x_0,\cdots, x_n]_{\lambda,g})\in {\mathbb H}_{n,d}^{\lambda, \Gamma} (P_d (\Gamma)),$$
 where the summation  $\sum_{(x_0,\cdots, x_n)}$ is taken over all ordered $n$-simplices of $ (P_d (\Gamma))_{g,d}$.
We remark that the choice of $n$ guarantees that the trace  in the above definition of the Connes-Chern character is finite.

By Propositions 4.4 and 4.5, the above Connes-Chern character induces a Connes-Chern character:

$$c_n:  K_0( {\cal S}_p\Gamma)\rightarrow
\lim_{d \rightarrow \infty}( \oplus _{k~even, ~k\leq n} ~{\mathbb H}_{k}^{ \Gamma} (P_d (\Gamma))).$$

For an arbitrary non-negative integer $n$, let $n'$ be a positive even integer satisfying $n'\geq max\{n, p\}.$
Let $\pi_{n', n}$ be the natural projection from \linebreak
$\lim_{d \rightarrow \infty}( \oplus _{k~even, ~k\leq n'} ~{\mathbb H}_{k}^{ \Gamma} ( P_d (\Gamma)))$
to $\lim_{d \rightarrow \infty}( \oplus _{k~even, ~k\leq n} ~{\mathbb H}_{k}^{ \Gamma} ( P_d (\Gamma)))$. We define the Connes-Chern character
$c_n$ from $ K_0( {\cal S}_p\Gamma)$ to \linebreak
$\lim_{d \rightarrow \infty}( \oplus _{k~even, ~k\leq n} ~{\mathbb H}_{k}^{ \Gamma} ( P_d (\Gamma)))$ to be $\pi_{n', n} \circ c_{n'}$. It is not difficult to verify that the definition of $c_n$ is independent of the choice of $n'$.

The proof of the following proposition is similar to the proof of Proposition 4.1 and is therefore omitted.

\begin{prop}
Let $\Gamma$ be a countable group. For any non-negative integer $n$,
the Connes-Chern character $c_n$ is a well defined homomorphism from $ K_0( {\cal S}_p\Gamma)$ to
$\lim_{d \rightarrow \infty}( \oplus _{k~even, ~k\leq n} ~{\mathbb H}_{k}^{ \Gamma} ( P_d (\Gamma)))$.
\end{prop}

Using the definition of lower algebraic K-theory, for each $p\geq 1$ and any non-negative integer $n$, we can similarly define
$$c_n: K_i( {\cal S}_p\Gamma) \rightarrow  \lim_{d\rightarrow \infty} (\oplus _{k+i~even, k\leq n~} ~{\mathbb H}_{k}^{ \Gamma} (  P_d (\Gamma) ))$$
for each $i<0$.

Finally, with the help of  the equality ${\cal S}=\cup _{p\geq 1}{\cal S}_p$, we obtain a Connes-Chern character
$$c_n: K_i( {\cal S}\Gamma) \rightarrow  \lim_{d\rightarrow \infty} (\oplus _{k+i~even, ~k\leq n}~ {\mathbb H}_{k}^{ \Gamma} ( P_d (\Gamma) ))$$
for each non-negative integer $n$ and  $i\leq 0$.

Notice that when $\Gamma$ is torsion free, $\Gamma_{fin}$ consists only of the identity element and the above definition coincides with the prior definition for the torsion free case.

In the rest of this section, we study a local property of  the Connes-Chern character for K-theory elements with small propagations.
This local property will play an important role in the proof of the main theorem of this paper.

We shall need a few preparations to explain the concept of propagation in a continuous setting.
Let $X$ be a $\Gamma$-invariant simplicial subspace of $P_{d_0} (\Gamma) $ for some $d_0\geq 0$. Endow $P_{d_0}(\Gamma)$ with a metric $d$ such that its restriction to each simplex is the standard metric and $d(\gamma_1, \gamma_2)\leq d_{\Gamma} (\gamma_1, \gamma_2)$ for all $\gamma_1$ and $\gamma_2$ in $\Gamma \subseteq P_{d_0}(\Gamma)$, where $d_{\Gamma}$ is the proper length metric on $\Gamma$. Let  $X$ be given     the  simplicial  metric of $P_{d_0} (\Gamma)$.
Let $H$ be a Hilbert space with a $\Gamma$-action and let $\phi$ be a $\ast$-homomorphism  from $C_0(X)$ to $B(H)$ which  is covariant in the sense that $\phi(\gamma f) h= (\gamma(\phi(f)) \gamma^{-1})h$ for all $\gamma\in \Gamma$, $f\in C_0(X)$ and $h\in H$.  Such a triple $(C_0(X), \Gamma, \phi)$ is called a covariant system.

The following definition is due to John Roe [Roe].

\begin{defi}  Let $H$ be a Hilbert space  and let $\phi$ be a $\ast$-homomorphism  from $C_0(X)$ to $B(H),$
the $C^*$-algebra of all bounded operators on $H$.
Let $T$ be a bounded linear operator acting on $H$.
\begin{enumerate}
\item[(1)]
The support of $T$ is defined to be the complement (in $X\times X$) of the set of all points
$(x,y) \in X\times X$  for which there exists $f\in C_0(X)$ and $g\in C_0(X)$ satisfying
$\phi(f)T\phi(g)=0$ and $f(x)\neq 0$ and $g(y)\neq 0$;
\item[(2)] The propagation of $T$ is defined to be:
$$\sup ~ \{ d(x,y): (x, y)\in {\rm Supp}(T)\};$$
\item[(3)]  Given $p\geq 1$, $T$ is said to be locally Schatten $p$-class if $\phi(f)T$ and $T\phi(f)$ are Schatten $p$-class operators for each $f\in C_c(X)$, the algebra of all compactly supported continuous functions on $X$.
\end{enumerate}
\end{defi}

\begin{defi}
We define  the covariant system $(C_0(X), \Gamma, \phi)$ to be admissible if
\begin{enumerate}
\item[(1)] the $\Gamma$-action on $X$ is proper and cocompact;
\item[(2)] $\phi$ is nondegenerate in the sense that $\phi( C_0(X))H$ is dense in $H$;
\item[(3)] $\phi(f)$ is noncompact for any nonzero function $f\in C_0(X)$;
\item[(4)] for each $x\in X$,  the action of  the stabilizer group $\Gamma_x$ on $H$ is regular in the sense that
it is isomorphic to the action of $\Gamma_x$ on $l^2(\Gamma_x)\otimes  W$ for some infinite dimensional  Hilbert space $W$, where the $\Gamma_x$ action on $l^2(\Gamma_x)$ is regular and the $\Gamma_x$ action on $W$  is trivial.
\end{enumerate}
\end{defi}

We remark  that condition (4) in the above definition is unnecessary if $\Gamma$ acts on $X$ freely.
In particular, if $M$ is a compact manifold and $\Gamma =\pi_1(M)$, then $(C_0(\widetilde{M}), \Gamma, \phi)$ is an admissible covariant system, where $\widetilde{M}$ is the universal cover of $M$ and $\phi(f)\xi=f\xi$ for each $f\in C_0(\widetilde{M})$ and all $\xi\in L^2(\widetilde{M}).$ In general, for each locally compact metric space with a proper and cocompact isometric action of $\Gamma$, there exists an admissible covariant system $(C_0(X), \Gamma, \phi)$.

\begin{defi}
For any $p\geq 1$, let $(C_0(X), \Gamma, \phi)$ be an admissible covariant system. We define $ \mathbb{C}_p(\Gamma, X, H)$
to be the  ring of  $\Gamma$-invariant locally Schatten $p$-class operators
 acting on $H$ with finite propagation.
\end{defi}

The proof of the following useful result is straightforward and is therefore omitted.

\begin{prop} Let $\Gamma$ be a countable group. Let $X$ be a simplicial complex with a simplicial proper and cocompact action of $\Gamma$.
If $(C_0(X), \Gamma, \phi)$ is an admissible covariant system, then the ring $\mathbb{C}_p(\Gamma, X,  H)$ is isomorphic to the ring
${\cal S}_p\Gamma$.
\end{prop}

For each $r>0$, let $X_{r}$ be a  $\Gamma$-invariant discrete subset of $X$ such that

\noindent (1) $X_{r}$ has bounded geometry, i.e. for each $R>0$, there exists $N>0$ such that any ball in $X_{r}$ with radius $R$ has at most $N$ elements;

\noindent (2) $X_{r}$ is
$r$-dense in $X$, i.e. $d(x, X_{r})<r$ for every $x\in X$;

\noindent (3) $X_{r}$ is uniformly discrete, i.e. there exists $k_r>0$ such that $d(z,z')\geq k_r$ for all distinct pairs of elements
$z$ and $z'$ in $X_{r}$.

 Let $\{U_z\}_{z\in X_{r}} $ be a $\Gamma$-equivariant disjoint Borel cover of $X$ such that $z\in U_z$ and
 $diameter(U_z)<r$ for all $z$.
Let $\chi_z$ be the characteristic function of $U_z$. Extend the $\ast$-representation $\phi$ to the algebra of all bounded Borel functions.
If $k\in \mathbb{C}_p(\Gamma, X,  H)$, let $k(x,y) =\phi (\chi_x) k \phi (\chi_y)$ for all $x$ and $y$ in $X_{r}.$

For any $r>0$, let $U_z'$ be  the $10r$-neighborhood of $U_z$ for each $z\in X_{r}$, i.e.
$$U_z'=\{ x\in X: d(x, U_z)<10 r\}.$$
Let $O_{r}(X)=\{U_z'\}_{z\in X_{r}} $.  Note that $O_{r}(X)$ is an open cover of $X$.

Let $ N(O_{r}(X))$ be the nerve space of the open cover $O_{r}(X)$.  We equip the vertex set $V$ of the simplicial complex $N(O_r (X))$  with the pseudo metric $d_V$ defined by:
$$ d_V(W, W')= \sup\{d(x,y): x\in W, y\in W'\} $$ for any pair of vertices $W$ and $W'$ in $N(O_r (X))$.

For each non-negative even integer $n\geq p$, we  define the Connes-Chern character
$$c_n: K_0( \mathbb{C}_p(\Gamma, X,  H)) \rightarrow \oplus _{k~even, k\leq n} ~{\mathbb H}_{k,r}^{ \Gamma} ( N(O_{r}(X)) )$$
as follows.

Let $\tilde{q}$ be an idempotent  in $M_m(\mathbb{C}_p(\Gamma, X,  H) )^+)$ and $\tilde{q}=q+q_0$ for some $q\in M_m(\mathbb{C}_p(\Gamma, X,  H))$ and    idempotent $q_0\in M_m ({\mathbb C})$. Let $r$ be the propagation of $q$.
 Let $n$ be an even integer satisfying $n\geq p$.

The Connes-Chern character  of
$[\tilde{q}] -[q_0]$ is defined to be homology class of
$$\sum_{g\in \Gamma_{fin}} (g, \sum_{(x_0,\cdots, x_n) } tr(q(x_0, x_1) \cdots q(x_n, g^{-1}x_0))
[x_0,\cdots, x_n]_{\lambda,g})$$ in the homology group $ {\mathbb H}_{n,(n+1)r}^{\lambda,  \Gamma} ( N(O_{(n+1)r}(X)) )$,
where,  for each $g$ and $r$,
$(x_0,\cdots, x_n)$ denotes the ordered simplex $(U'_{x_0}, \cdots, U'_{x_n})$ in the  space \linebreak $(N(O_{(n+1)r}(X)))_{g,r}$,
$[x_0,\cdots, x_n]_{\lambda,g}$ denotes
the $g$-cyclically oriented simplex $[U'_{x_0}, \cdots, U'_{x_n}]_{\lambda, g}$  in
$ (N(O_{(n+1)r}(X)))_{g,r} $, and
the summation $\sum_{(x_0,\cdots, x_n) }$ in the above formula
is taken over all ordered simplices in $ (N(O_{(n+1)r}(X)))_{g,r} $.

The following proposition follows from the above definition and the proof of Proposition 4.1.

\begin{prop} Let $\Gamma$ be a countable group. Let $X$ be a simplicial complex with a simplicial proper and cocompact action of $\Gamma$.
 For each $r>0$, let $n$ be an even integer satisfying $n\geq p$.
If $(C_0(X), \Gamma, \phi)$ is an admissible covariant system, then
the Connes-Chern character $c_n$ of an element in $ K_0(\mathbb{C}_p(\Gamma, X,  H) )$ with propagation less than or equal to  $r>0$ is a homology class in
${\mathbb H}_{n,r}^{ \lambda, \Gamma} ( N(O_{(n+1)r}(X)) )$.
\end{prop}

We identify $K_0( {\cal S}_p\Gamma)$  with  $\lim_{d\rightarrow \infty}\lim_{X} K_0( \mathbb{C}_p(\Gamma, X,  H))$ using Proposition 4.10,
where the direct limit  $\lim_{X}$ is  taken over the directed system of all $\Gamma$-invariant, $\Gamma$-compact
subsets of $P_d(\Gamma)$. We also identify  \linebreak
$\lim_{d\rightarrow \infty} \lim_{X}  ~{\mathbb H}_{n,r}^{\lambda,  \Gamma} (N(O_{r}(X)) ))$ with
$\lim_{ d \rightarrow \infty}( \oplus_{k~even, k\leq n}
~{\mathbb H}_k^{ \Gamma} (P_d(\Gamma) ))$ using \linebreak Propositions 4.4 and 4.5,  where the direct limit  $\lim_{X}$ is again  taken over the directed system of all $\Gamma$-invariant, $\Gamma$-compact
subsets of $P_d(\Gamma)$.
 Using the projection $\pi_{n',n}$ from $\lim_{ d \rightarrow \infty}( \oplus_{k~even, k\leq n'}
~{\mathbb H}_k^{ \Gamma} (P_d(\Gamma) ))$ to $\lim_{ d \rightarrow \infty}( \oplus_{k~even, k\leq n}
~{\mathbb H}_k^{ \Gamma} (P_d(\Gamma) ))$ for  $n'=max\{n,p\}$, the above Connes-Chern character induces a Connes-Chern character:
$$c_n=\pi_{n',n} \circ c_{n'}: K_0( {\cal S}_p\Gamma) \rightarrow  \lim_{d\rightarrow \infty} (\oplus _{k~even,~k\leq n}~~{\mathbb H}_k^{ \Gamma} (P_d (\Gamma) ))$$
for any non-negative even integer $n$.
This construction gives back the  Connes-Chern character in Proposition 4.6.

In the following corollary, we demonstrate a local property of the Connes-Chern character.
This local property of the Connes-Chern character  plays an important role in the proof of Theorem 1.1.

\begin{cor} Let $\Gamma$ be a countable group. Let $X$ be a simplicial complex with a simplicial proper and cocompact action of $\Gamma$.
 Let $\tilde{q}=q+q_0$ be an element in $ M_m( \mathbb{C}_p(\Gamma, X,  H)^+)$ such that $q\in  M_m (\mathbb{C}_p(\Gamma, X,  H))$, $q_0\in M_m({\mathbb C})$, $\tilde{q}$ and $q_0$ are idempotents. For any non-negative integer $n$,
when  the propagation of $q$ is sufficiently small, $c_n([\tilde{q}]-[q_0])$ can be represented by a homology class in
$\oplus _{k~even,~k\leq n} ~{\mathbb H}_k^{ \Gamma} (X)$. More generally for each $i\leq 0$ and $p\geq 1$,
the Connes-Chern character $c_n$ of an element in $K_i (\mathbb{C}_p(\Gamma, X,  H)) $ with sufficiently small propagation can be represented by a homology class in
$\oplus _{k~even,~ k\leq n}~ {\mathbb H}_{k+i}^{ \Gamma} (X)$.
\end{cor}
\proof
Let $\tilde{q}=q+q_0$ be an element in $ M_m( \mathbb{C}_p(\Gamma, X,  H)^+)$ such that $q\in  M_m (\mathbb{C}_p(\Gamma, X,  H))$, $q_0\in M_m({\mathbb C})$, $\tilde{q}$ and $q_0$ are idempotents.
If $n\geq p$ and the propagation of $q$ is less than or equal to  $r>0$, then by Proposition 4.11 we know that $c_n([\tilde{q}]-[q_0])$ can be represented by a homology class in
$ {\mathbb H}_{n,r}^{\lambda, \Gamma} (N(O_{(n+1) r}(X)) )$.

Let
$$r_0 =\inf ~\{ ~d(gU'_z, U_z'): ~z\in X_\epsilon, g\in \Gamma, gU_z'\neq U_z'\},$$
where   $\{U'_z\}_{z\in X_r}$ is the open cover $O_r(X)$ in the definition of the Connes-Chern character.

Assume that   $r$ is  a sufficiently small positive number for the rest of this proof. We have $r_0>0$. We  choose $\{U_z\}_{z\in X_r}$
 such that $\{U_z'\}_{z\in X_r}$  is a good cover.

If $r<r_0$,  then we have the following:

\noindent{[1]} by the definition of $ {\mathbb H}_{n,r}^{\lambda, \Gamma} (N(O_{r}(X)) )$, the homology group $ {\mathbb H}_{n,r}^{\lambda,  \Gamma} ( N(O_{r}(X)) )$ is equal to $H_n^{\lambda, \Gamma} (\widehat{( N(O_{r}(X)))})$;

\noindent{[2]} by Proposition 4.2, the homology group $H_n^{\lambda, \Gamma} (\widehat{( N(O_{r}(X)))})$
can be identified with
 $ \oplus _{k~even, k\leq n}~  H_{k}^{ \Gamma} ( \widehat{(N(O_{r}(X))} ); $

\noindent{[3]} by the choice of  $\{U'_z\}_{z\in X_r}$, the homology group $H_{k}^{ \Gamma} ( \widehat{(N(O_{r}(X)))} )$ is equal to
  ${\mathbb H}_k^{ \Gamma} (X)$ for each $k$.

When $n\geq p$, Corollary 4.12  follows from the above statements and the definition of the  lower algebraic K-groups. The case of an arbitrary
non-negative integer $n$ can be reduced to this special case by considering $n'=max\{n, p\}$ and using the identity $c_n =\pi_{n', n}\circ c_{n'}$, where
$\pi_{n', n}$ is the projection from $\oplus _{k~even,~ k\leq n'}~ {\mathbb H}_{k+i}^{ \Gamma} (X)$ to $\oplus _{k~even,~ k\leq n}~ {\mathbb H}_{k+i}^{ \Gamma} (X)$.
\qed

\section{Proof of the main result}

In this section, we give a proof of Theorem 1.1.

By Proposition 2.4, Theorem 1.1 follows from the following result.

\begin{theo} Let ${\cal S}$ be the ring of all  Schatten class operators on an infinite dimensional and separable Hilbert space.  The assembly map
$$A: H_n^{Or\Gamma}(E_{\cal FIN} (\Gamma), {\mathbb K}({\cal S})^{-\infty}) \rightarrow K_n ({\cal S}\Gamma)$$ is rationally injective
for any  group $\Gamma$.
\end{theo}

\proof Without loss of generality, we can assume that $\Gamma$ is countable (this is because every group is an inductive limit of countable groups).
 We recall that $E_{\cal FIN} (\Gamma)$ can be identified with $\cup_{d\geq 1} P_d(\Gamma)$.
For each $i\leq 0 $ and non-negative integer $n$, composing the assembly map
$$A: H_i^{Or\Gamma}(E_{\cal FIN} (\Gamma), {\mathbb K}({\cal S})^{-\infty}) \rightarrow K_i ({\cal S}\Gamma)$$
with the Connes-Chern character
$$c_n:~~ K_i( {\cal S}\Gamma) \longrightarrow \lim_{d\rightarrow \infty}( \oplus _{k~even, ~k\leq n} ~ {\mathbb H}_{k+i}^{ \Gamma} ( P_d (\Gamma))),$$
we obtain a homomorphism
$$\psi_{i,n}:~~ H_i^{Or\Gamma}(E_{\cal FIN} (\Gamma), {\mathbb K}({\cal S})^{-\infty}) \rightarrow  \lim_{d\rightarrow \infty}( \oplus _{k~even, k\leq n}~  {\mathbb H}_{k+i}^{\Gamma} ( P_d (\Gamma))).$$

By using the fact that $H_i^{Or\Gamma}(E_{\cal FIN} (\Gamma), {\mathbb K}({\cal S})^{-\infty})$ is \linebreak
$\lim_{d\rightarrow \infty}H_i^{Or\Gamma}(P_{d}(\Gamma), {\mathbb K}({\cal S})^{-\infty})$, we obtain a homomorphism
$$\psi_{i}:~~ H_i^{Or\Gamma}(E_{\cal FIN} (\Gamma), {\mathbb K}({\cal S})^{-\infty}) \rightarrow  \lim_{d\rightarrow \infty}( \oplus _{k~even}~  {\mathbb H}_{k+i}^{\Gamma} ( P_d (\Gamma)))$$ such that
$\psi_i$ coincides with $\psi_{i, n}$ on the image of the natural map from  \linebreak
$H_i^{Or\Gamma}(P_d(\Gamma), {\mathbb K}({\cal S})^{-\infty})$ to
$H_i^{Or\Gamma}(E_{\cal FIN} (\Gamma), {\mathbb K}({\cal S})^{-\infty})$ when $n+i$ is greater than or equal to the dimension of $P_d(\Gamma)$.
Such a map $\psi_i$ is unique and is independent of the choice of $n$.

For each simplicial  $\Gamma$-invariant and $\Gamma$-cocompact
subspace $X$ of $E_{\cal FIN} (\Gamma)$,
let $\mathbb{C}(\Gamma, X,  H) =\cup_{p\geq 1} \mathbb{C}_p(\Gamma, X,  H),$ where $\mathbb{C}_p(\Gamma, X,  H)$ is as in definition 4.9.
By the definition of the assembly map in [BFJR]  and the fact that this assembly map coincides with the classic assembly map (Corollary 6.3 in [BFJR]),  K-theory elements in the image of the assembly map
$$A: H_i^{Or\Gamma}(X, {\mathbb K}({\cal S})^{-\infty}) \rightarrow K_i (\mathbb{C}(\Gamma, X,  H))$$
can be represented by elements with arbitrarily small propagation for $i\leq 0$.

This, together with Corollary 4.12, implies that
there exists a map (still denoted by $\psi_i$)
$$\psi_i:~~ H_i^{Or\Gamma}(X, {\mathbb K}({\cal S})^{-\infty}) \rightarrow  ( \oplus _{k~even}~  {\mathbb H}_{k+i}^{ \Gamma} ( X))$$ for each non-positive integer $i$ such  that the following diagram commutes:
$$\begin{array}[c]{ccc}
H_i^{Or\Gamma}(X, {\mathbb K}({\cal S})^{-\infty})  &\stackrel{\psi_i}{\longrightarrow}&  ( \oplus _{k~even} ~ {\mathbb H}_{k+i}^{\Gamma} ( X)) \\
\downarrow\scriptstyle{j_\ast}&&\downarrow\scriptstyle{j'_\ast}\\
H_i^{Or\Gamma}(E_{\cal FIN} (\Gamma), {\mathbb K}({\cal S})^{-\infty}) &\stackrel{\psi_i}{\longrightarrow}&  \lim_{d\rightarrow \infty}( \oplus _{k~even} ~ {\mathbb H}_{k+i}^{\Gamma} ( P_d (\Gamma)))
\end{array},$$ where $j_\ast$ and $j'_\ast$ are respectively induced by the inclusion maps.

If $X=\Gamma/F$ as $\Gamma$-spaces for some finite subgroup $F$ of $\Gamma$, then it is straightforward to verify that $\psi_i$ is an isomorphism after tensoring with ${\mathbb C}$. In fact, both sides are naturally isomorphic to the group $R(F)\otimes {\mathbb C}$, where $R(F)$ is the representation ring of $F$ viewed as an additive group.
Recall that, by Proposition 7.2.3, Remark 7.2.6 and  Theorem 8.2.5 in [CT], we have $K_n ({\cal S}) ={\mathbb Z}$ when $n$ is even and $K_n ({\cal S}) =0$ when $n$ is odd.
As a consequence, the homology theory $H_i^{Or\Gamma}(X, {\mathbb K}({\cal S})^{-\infty})$ is 2-periodic.
 Note that the homology theory $\oplus _{k~even} ~ {\mathbb H}_{k+i}^{ \Gamma} (X)$ is 2-periodic by definition.
By the proof of the Mayer-Vietoris sequence using the mapping cone and the definition of the Connes-Chern character, we know the homomorphisms  $\psi_i$ commute with the Mayer-Vietoris sequences (up to scalars)

Using the above results, the fact that both homology theories satisfy the Mayer-Vietoris sequence  and  a five lemma argument, we can prove that
the map $$\psi_i:~~ H_i^{Or\Gamma}(X, {\mathbb K}({\cal S})^{-\infty}) \rightarrow  ( \oplus _{k~even} ~ {\mathbb H}_{k+i}^{ \Gamma} ( X))$$
  is  an isomorphism  after tensoring with ${\mathbb C}$ for $i\leq 0$. This implies that the assembly map $A$ is rationally injective for $i\leq 0$.
Now Theorem 5.1 follows from Proposition 3.1.
\qed

We  comment that the algebraic K-theory isomorphism conjecture for the ring ${\cal S} \Gamma$ can be viewed as an algebraic counterpart of the Baum-Connes conjecture  for  the K-theory of the reduced group $C^\ast$-algebra of $\Gamma$ [BC1]. The Farrell-Jones isomorphism conjecture and the Baum-Connes conjecture imply the following conjecture.

 \begin{conj} Let $K$ be  the $C^\ast$-algebra of all compact operators on an infinite dimensional and separable  Hilbert space, let   $C_r^\ast (\Gamma)$ be the reduced group $C^\ast$-algebra of $\Gamma$.   The natural homomorphism  $$i_\ast: K_n (  {\cal S} \Gamma)\rightarrow K_n (C_r^\ast (\Gamma)\otimes K)$$
 is an isomorphism,
 where $C_r^\ast (\Gamma)\otimes K $ is the $C^\ast$-algebraic tensor product of   $C_r^\ast (\Gamma)$ with $K$
 and $i$ is the inclusion map from $ {\cal S} \Gamma$ to $C_r^\ast (\Gamma)\otimes K $.
 \end{conj}

  We remark that, by a theorem of Suslin-Wodzicki [SW],  the algebraic K-theory $K_n (C_r^\ast (\Gamma)\otimes K)$ is isomorphic to the topological K-theory $K_n ^{top}(C_r^\ast (\Gamma))$.
  Theorem 1.1 implies that the Novikov higher signature conjecture follows from the (rational) injectivity  of $i_\ast$ in the above conjecture.

   Finally we  speculate that the (algebraic) bivariant K-theory of Cuntz, Cuntz-Thom and Corti\~nas-Thom should be useful in studying the algebraic K-theory isomorphism conjecture for ${\cal S}\Gamma$ [Cu1] [Cu2] [CuT] [CT1].

\noindent Shanghai Center for Mathematical Sciences, China.

\noindent Department of Mathematics, Texas A\&M University,

\noindent College Station, TX 77843, USA.

\noindent e-mail: guoliangyu@math.tamu.edu

\end{document}